\documentclass[a4paper,twocolumn]{article} 

\usepackage{graphics} 
\usepackage{epsfig} 
\usepackage{times} 
\usepackage{amsmath} 
\usepackage{amssymb}  
\usepackage[left=2cm, right=2cm, top=3cm, bottom=3cm]{geometry}

\usepackage{epstopdf}
\usepackage{theorem}
\usepackage[]{units}
\usepackage{subfig}
\usepackage{psfrag}
\usepackage{graphicx}
\usepackage{mathabx}
\usepackage{ stmaryrd }
\usepackage{wrapfig}
\usepackage{color}
\usepackage[percent]{overpic}

\usepackage{multirow}

 \newtheorem{lemma}{Lemma}
 \newtheorem{corollary}{Corollary}

  \newtheorem{proof}{Proof}
 
 \newcommand{\R}{\mathbb{R}}
\newcommand{\N}{\mathbb{N}}

\title{Reducing the computational effort of min-max model predictive control with regional feedback laws
}

\author{Kai K\"onig and Martin M\"onnigmann\thanks{Corresponding author.}\\
Automatic Control and Systems Theory, Department of Mechanical Engineering,\\
	Ruhr-Universit\"at Bochum, 44801 Bochum, Germany.\\ E-mail: {\tt\small kai.koenig-h4d@rub.de} and {\tt\small martin.moennigmann@rub.de}}

\begin{document}

\maketitle

\begin{abstract}
Recently, a regional MPC approach has been proposed that exploits the piecewise affine structure of the optimal solution (without computing the entire explicit solution before).
Here, {\em{regional}} refers to the idea of using the affine feedback law that is optimal in a vicinity of the current state of operation, and therefore provides the optimal input signal without requiring to solve a QP. 
In the present paper, we apply the idea of regional MPC to min-max MPC problems. 
We show that the new robust approach can significantly reduce the number of QPs to be solved within min-max MPC resulting in a reduced overall computational effort. 
Moreover, we compare the performance of the new approach to an existing robust regional MPC approach using a numerical example with varying horizon. 
Finally, we provide a rule for choosing a suitable robust regional MPC approach based on the horizon.
\end{abstract}

\section{Introduction}
Model predictive control (MPC) is an established method for the control of dynamical systems with state and input constraints. The method, however, requires the solution of an optimization problem in every time step resulting in a high computational effort. 
The effort is even higher if disturbances have to be taken into account explicitly. One approach to cope with disturbances is to minimize the cost function for the worst possible case of disturbances, which is known as min-max MPC (see, e.g., \cite{Witsenhausen1968, Campo1987,Scokaert1998, Bemporad2003, Kerrigan2004}). Some approaches propose to minimize over a sequence of control corrections to a given linear feedback law for the nominal plant (see, e.g., \cite{Bemporad1998, Lofberg2003, Munoz2007}). By this, the conservatism of the open-loop predictions is reduced without increasing the computational effort.  
 
Recently, a regional MPC approach has been proposed in which an optimal affine feedback law and its polytopic region of validity are computed from the solution at the current state and reused for subsequent states to avoid online optimizations \cite{Jost2015a}.  
Only when leaving the current polytope, a new optimization is executed. To reduce the number of optimizations further, several extensions, such as using active set updates and suboptimality, have been developed (see, e.g., \cite{Koenig2017a, Koenig2017d, Koenig2017c}).  
Since the original approach from \cite{Jost2015a} does not take disturbances into account directly (note it is robust to some extent due to the polytopes), in \cite{SchulzeDarup2017CDC} the underlying classical MPC problem has been replaced by a tube-based MPC problem (see, e.g., \cite{MAYNE2005}) resulting in a robust regional MPC approach. 
The robust approach has also been extended by active set updates and suboptimality and it has been shown that the number of online optimizations in tube-based MPC can be reduced significantly \cite{BernerP2019b}.
Regional MPC has been implemented in a networked MPC variant on standard industrial hardware \cite{Berner2020} and even low-cost embedded hardware \cite{Koenig2020}.  

We apply the regional idea from \cite{Jost2015a} to a min-max optimization problem as proposed in \cite{Munoz2007}. This results in a new robust regional MPC approach, which we also extend by active set updates and suboptimality to reduce the number of online optimizations. We investigate the performance of the new approach using a numerical example with varying horizon and compare it to the tube-based approach from \cite{SchulzeDarup2017CDC} and the original approach from \cite{Jost2015a}. 
We show that the number of online optimizations and thus the overall computational effort in min-max MPC can be reduced significantly.
Moreover, we provide a horizon-dependent rule for choosing a suitable robust regional MPC approach, when the inherent robustness of the regional approach from \cite{Jost2015a} is not sufficient and robustness has to be guaranteed.

We state the system and problem class along with some preliminaries in Section \ref{sec:Problem}. The new approach is presented in Section \ref{sec:RMPC} and applied to an example in Section \ref{sec:Example}. In Section \ref{sec:Conclusion} we give a short summary.  

\subsection*{Notation}
For an arbitrary matrix $M \in \R^{a \times b}$, let $M^\mathcal{L}$ with $\mathcal{L} \subseteq \{1, \ldots, a\}$ describe the submatrix with the rows indicated by $\mathcal{L}$. If $\mathcal{L}$ contains only one element, say $\mathcal{L}= \left\{i\right\}$, we write $M^i$. A \textit{polytope} is a bounded set defined by a finite number of intersecting halfspaces $\mathcal{P}=\{x \in \R^n \ | \ Tx \leq d \}$ with $T \in \R^{r \times n}$ and $d \in \R^{r}$.

\section{Problem statement}\label{sec:Problem}

Consider a linear discrete-time system with bounded additive disturbances 
\begin{align}\label{eq:Sys}
\begin{split}
x(k+1)&=Ax(k)+Bu(k)+Dw(k),\, \quad x(0) \text{ given},
\end{split}
\end{align}
with states $x(k) \in \R^n$, inputs $u(k) \in \R^m$, disturbances $w(k) \in \R^s$ and system matrices $A \in \R^{n \times n}$, $B \in \R^{n \times m}$ and $D \in \R^{n \times s}$. Assume state and input constraints and bounded disturbances 
\begin{align}\label{eq:constraints}
x(k) \in \mathcal{X}, \quad u(k) \in \mathcal{U}, \quad w(k) \in \mathcal{D}
\end{align} 
apply for all $k \geq 0$, where $\mathcal{X}$, $\mathcal{U}$ and $\mathcal{D}$ are compact polyhedra that contain the origin as an interior point. We choose the control input in \eqref{eq:Sys} to be 
\begin{align}\label{eq:LQR_correction}
u(k)=-K_{\infty}x(k)+v(k)
\end{align}
with the LQR feedback gain $K_{\infty}$ and the correction control input $v(k) \in \R^m$ as in \cite{Munoz2007}. The dynamics of the system \eqref{eq:Sys} can be rewritten as
\begin{align}\label{eq:closedLoopSystem}
\begin{split}
x(k+1)&=A_{\mathrm{cl}} x(k)+Bv(k)+Dw(k) \quad \text{with}\\
A_{\mathrm{cl}}&=(A-BK_{\infty}).
\end{split}
\end{align}
The proposed MPC controller aims at steering system \eqref{eq:closedLoopSystem} to a robust positively invariant set (RPI set) around the origin while satisfying constraints \eqref{eq:constraints}. For this purpose, a min-max MPC problem is solved, in which a cost function is minimized for the worst possible case that may result from the disturbances. 
More precisely, the min-max problem 
\begin{align}
\small
\begin{split}
\min \limits_{X,V} \ &\max \limits_{W \in \mathcal{D}_{\mathrm{N}}} \ \tilde{x}(N)^{\prime}P\tilde{x}(N)+\sum \limits_{i=0}^{N-1} (\tilde{x}(i)^{\prime}Q\tilde{x}(i)+\tilde{u}(i)^{\prime}R\tilde{u}(i)) \\
\text{s.t.} \quad &\tilde{x}(0)=x,\\
&\tilde{x}(i+1)=A_{\mathrm{cl}}\tilde{x}(i)+B\tilde{v}(i)+D\tilde{w}(i),  \quad i=0, \ldots, N-1,\\
&\tilde{u}(i)=-K_{\infty}\tilde{x}(i)+\tilde{v}(i),  \quad \quad \quad \quad i=0, \ldots, N-1,\\
&\tilde{x}(i) \in \mathcal{X}, \quad \tilde{u}(i) \in \mathcal{U}, \quad \tilde{w}(i) \in \mathcal{D},  \quad i=0, \ldots, N-1,\\
&\tilde{x}(N) \in \mathcal{T}
\end{split}
\label{eq:MPCProblem}
\end{align} 
is solved for the current state $x$ on a receding horizon $N \in \N$, where $W=(\tilde{w}(0)^{\prime},\tilde{w}(1)^{\prime}, \ldots, \tilde{w}(N-1)^{\prime})^{\prime}$ represents a sequence of disturbances to the system and $X=(\tilde{x}(1)^{\prime},\tilde{x}(2)^{\prime}, \ldots, \tilde{x}(N)^{\prime})^{\prime}$ and $V=(\tilde{v}(0)^{\prime},\tilde{v}(1)^{\prime}, \ldots, \tilde{v}(N-1)^{\prime})^{\prime}$ are a predicted state sequence and correction control input sequence, respectively. Let $\mathcal{X}_f$ refer to the set of initial states $x$ for which \eqref{eq:MPCProblem} has a solution. The set $\mathcal{D}_{\mathrm{N}}$ is the set of possible disturbance sequences of length $N$, i.e., 
{\small
\begin{align}
\mathcal{D}_{\mathrm{N}}=\{(w(0), \ldots, w(N-1)) \ | \ w(i) \in \mathcal{D} \ \forall i \in \{0, \ldots, N-1  \}  \}.
\end{align} 
}
Note that the predicted states and control inputs in the cost function depend on the disturbances and control corrections due to the equality constraints. A closed-loop system operates by solving the problem \eqref{eq:MPCProblem} at each sampling instant and applying the optimal correction control input $v^\star(0)$ to system \eqref{eq:closedLoopSystem}. Assume $Q \succeq 0$ and $R \succ 0$ with the obvious dimensions, $(A,B)$ is stabilizable and $(Q^{\frac{1}{2}},A)$ is detectable. Robust stability can be guaranteed by choosing $P$ as the solution of the discrete-time algebraic Riccati equation, which implies $P \succ 0$, and $\mathcal{T}$ to a robust control invariant set for the LQR regulated system. 
Since the cost function in \eqref{eq:MPCProblem} is convex in $W$, the maximum of the inner optimization for all possible uncertainties $W \in \mathcal{D}_{\mathrm{N}}$ is attained at least at one of the vertices $\text{Ver}(\mathcal{D}_{\mathrm{N}})$ \cite{Munoz2007}. 

By some straightforward reformulations (more details in \cite{Ramirez2006, Munoz2007}), OCP \eqref{eq:MPCProblem} can be expressed as a quadratic program
\begin{align}\label{eq:reformulatedProblem}
\begin{split}
\min \limits_{Z, \gamma} \ &\frac{1}{2}Z^{\prime}HZ+\gamma\\
\text{s.t.} \quad &G_{\mathrm{m}} Z + g_{\mathrm{m}} \gamma \leq W_{\mathrm{m}}+S_{\mathrm{m}} x, \\
&G_{\mathrm{c}} Z \leq W_{\mathrm{c}}+S_{\mathrm{c}} x
\end{split} 
\end{align}
with the optimizers $Z \in \R^{mN}$ and $\gamma \in \R$, where 
\begin{align}\label{eq:VandZ}
Z=V+H^{-1} L' x
\end{align} 
and $H \in \R^{mN \times mN}$, $G_{\mathrm{m}} \in \R^{r \times mN}$, $g_{\mathrm{m}} \in \R^{r}$, $W_{\mathrm{m}} \in \R^r$, $S_{\mathrm{m}} \in \R^{r \times n}$, $G_{\mathrm{c}} \in \R^{l \times mN}$, $W_{\mathrm{c}} \in \R^l$, $S_{\mathrm{c}} \in \R^{l \times n}$ and $L \in \R^{n \times mN}$. 
We refer to \cite{Ramirez2006} and \cite{Munoz2007} for more details about the computation of the matrices. Each of the first $r$ inequalities in \eqref{eq:reformulatedProblem} corresponds to one of the vertices $\text{Ver}(\mathcal{D}_N)$. At the optimum $(Z^\star, \gamma^\star)$ at least one of these $r$ constraints is active, i.e., fulfilled with equality. The second $l$ inequalities collect the state and input constraints. The total number of constraints depending on the horizon $N$ is analyzed in Section \ref{sec:NumberConstraints} in more detail.
Note that problem \eqref{eq:reformulatedProblem} is strictly convex because $H \succ 0$ and $\gamma$ is the maximum of a nonconstant linear function on a polytope~\cite{Ramirez2006, Munoz2007}. Thus, the optimizer $(Z^\star, \gamma^\star)$ is unique.

\section{A min-max regional MPC approach}\label{sec:RMPC} 
We propose a min-max regional MPC approach. Roughly speaking, the approach results from replacing the classical MPC problem used in the original approach from \cite{Jost2015a} with the robust problem \eqref{eq:MPCProblem} respectively \eqref{eq:reformulatedProblem}. As in \cite{Jost2015a}, we show that affine feedback laws can be computed from the solution at the current state. In contrast to \cite{Jost2015a}, the resulting feedback laws are guaranteed to be robust. We explain the algorithm at the end of this section. In Sections \ref{sec:AccActiveSetUp} and \ref{sec:AccSuboptimality} we extend the new approach by active set updates and suboptimality. We start by introducing some technical details in preparation of Lemma \ref{lem:robustRegionalFeedbackLaw}, which constitutes the basis of our new approach. Using the augmented decision variable
\begin{align}\label{eq:augmentedVariable} 
\epsilon=[Z^{\prime}, \ \gamma^{\prime}]^{\prime}
\end{align} 
problem \eqref{eq:reformulatedProblem} can be rewritten as       
\begin{align}\label{eq:rereformulatedMPCProblem}
\begin{split}
\min \limits_{\epsilon} \ &\frac{1}{2}\epsilon^{\prime}\hat{H}\epsilon+c^{\prime} \epsilon\\
\text{s.t.} \quad &G \epsilon \leq W+S x
\end{split} 
\end{align}
with the block matrices 
\begin{align}
\hat{H}=\begin{pmatrix} H & 0 \\ 0 & 0\end{pmatrix}, \quad c^{\prime}=\begin{pmatrix} 0 & 1 \end{pmatrix}
\end{align}
and $G \in \R^{(l+r) \times (mN+1)}$, $W \in \R^{(l+r)}$ and $S^{(l+r) \times n}$. Note that $\hat{H}$ is a singular matrix, because $\gamma$ enters only linearly (cf. \eqref{eq:reformulatedProblem}). Consequently, $\hat{H}$ is not invertible and a feedback law and its polytope must be computed with the null space method \cite{Nocedal2006, Tondel2003}. However, the optimizer $\epsilon^\star$ is uniquely defined due to uniqueness of the optimizer $(Z^\star, \gamma^\star)$ for \eqref{eq:reformulatedProblem}.  

We introduce the sets of active, inactive and weakly active constraints
\begin{align}\label{eq:sets}
\begin{split}
\mathcal{A}(x)&=\{i \in \mathcal{Q}~|~G^i\epsilon^\star(x)-W^i-S^ix=0 \},\\
\mathcal{I}(x)&=\{i \in \mathcal{Q}~|~G^i\epsilon^\star(x)-W^i-S^ix<0 \},\\
\mathcal{W}(x)&=\{i \in \mathcal{A}(x)~|~\lambda^{\star i}(x)=0 \}
\end{split}
\end{align}
with $\mathcal{I}(x)=\mathcal{Q}\backslash\mathcal{A}(x)$, where $\mathcal{Q}:=\{1, \ldots, r+l \}$ is the set of all constraint indices and $\lambda \in \R^{r+l}$ are the Lagrange multipliers. The following lemma is based on the results in \cite{Tondel2003}.

\begin{lemma}\label{lem:robustRegionalFeedbackLaw}
Let $x \in \mathcal{X}_{f}$ be arbitrary and $\mathcal{A}(x)$ the corresponding active set. Assume the matrix $G^{\mathcal{A}}$ has full row rank and let $(G^{\mathcal{A}})^{\prime}=\begin{bmatrix} E & J \end{bmatrix} \begin{bmatrix} F \\ 0\end{bmatrix}$ be a QR factorization of $(G^{\mathcal{A}})^{\prime}$. Let 
\begin{align}\label{eq:Kb}
\begin{split}
K_{\epsilon}&=E \Theta S^{\mathcal{A}}- J \Psi J^{\prime} \hat{H} E \Theta S^{\mathcal{A}},\\
b_{\epsilon}&=(E-J \Psi J^{\prime} \hat{H} E) \Theta W^{\mathcal{A}} - J \Psi J^{\prime} c, \\
K_{\lambda}&=-\Theta^{\prime} E^{\prime} (\hat{H} K_{\epsilon}), \\
b_{\lambda}&=-\Theta^{\prime} E^{\prime} (\hat{H} b_{\epsilon}+c),
\end{split} 
\end{align}
where $\Psi=(J^{\prime}\hat{H}J)^{-1}$ and $\Theta=(G^{\mathcal{A}}E)^{-1}$ exist by construction.
Then $\epsilon^\star(x)=K_{\epsilon}x+b_{\epsilon}$ is the affine optimizer and $\lambda^{\star \mathcal{A}}(x)=K_{\lambda}x+b_{\lambda}$ are the corresponding affine active Lagrange multipliers on the polytope $\mathcal{P}=\{x \in \R^n \ | \ Tx \leq d\}$ with
\begin{align}\label{eq:Td}
\begin{split}
T=\begin{pmatrix} G^{\mathcal{I}} K_\epsilon-S^{\mathcal{I}} \\ -K_{\lambda}\end{pmatrix} \quad \text{and} \quad
d=\begin{pmatrix} W^{\mathcal{I}} -G^{\mathcal{I}} b_{\epsilon} \\ b_{\lambda}\end{pmatrix}.
\end{split} 
\end{align}
\end{lemma}

\begin{proof}
Lemma \ref{lem:robustRegionalFeedbackLaw} can be proven analogously to Theorem 1 in \cite{Tondel2003}. We carry out the proof for completeness. The Karush-Kuhn-Tucker conditions for Problem \eqref{eq:rereformulatedMPCProblem} read 
\begin{align}
\hat{H} \epsilon +c+G' \lambda&=0, \label{eq:KKT_a_minmax} \\
\lambda^i (G^i \epsilon - W^i - S^i x)&=0, \quad i=1,\ldots,l+r, \label{eq:KKT_b_minmax}\\
G^i \epsilon-W^i-S^i x &\leq 0,\quad i=1,\ldots,l+r,\label{eq:KKT_c_minmax} \\
\lambda^i &\geq 0,  \quad i=1,\ldots,l+r.\label{eq:KKT_d_minmax}  
\end{align}
We devide the vector $\epsilon$ into two parts 
\begin{align}\label{eq:epsilon}
\epsilon=E \epsilon_{\mathrm{E}} + J \epsilon_{\mathrm{J}}.
\end{align}
Substituting \eqref{eq:epsilon} into those rows of \eqref{eq:KKT_c_minmax} that are fulfilled with equality, i.e., 
\begin{align}
G^{\mathcal{A}} \epsilon - W^{\mathcal{A}} - S^{\mathcal{A}} x=0,
\end{align}
yields
\begin{align}\label{eq:ConstraintsWithEpsilon}
G^{\mathcal{A}} E \epsilon_{\mathrm{E}}+ G^{\mathcal{A}} J \epsilon_{\mathrm{J}}-W^{\mathcal{A}} - S^{\mathcal{A}} x=0.
\end{align}
Since $G^{\mathcal{A}} J=0$ and $G^{\mathcal{A}} E$ is invertible by construction, relation \eqref{eq:ConstraintsWithEpsilon} results in 
\begin{align}\label{eq:epsilonY}
\epsilon_{\mathrm{E}}=(G^{\mathcal{A}} E)^{-1}(W^{\mathcal{A}} + S^{\mathcal{A}} x).
\end{align}
Substituting \eqref{eq:epsilon} into the stationarity condition \eqref{eq:KKT_a_minmax} and multiplying by $J^{\prime}$ results in
\begin{align}\label{eq:StationarityWithEpsilon}
J^{\prime} \hat{H} E \epsilon_{\mathrm{E}}+J^{\prime} \hat{H} J \epsilon_{\mathrm{J}} + J^{\prime} c + J^{\prime} (G^{\mathcal{A}})^{\prime} \lambda^{\mathcal{A}}  = 0,
\end{align}
where $G^{\prime} \lambda = (G^{\mathcal{A}})^{\prime} \lambda^{\mathcal{A}}$ holds because $\lambda^{\mathcal{I}}=0$. Solving \eqref{eq:StationarityWithEpsilon} for $\epsilon_{\mathrm{J}}$ with $J^{\prime} (G^{\mathcal{A}})^{\prime}=0$ by construction results in
\begin{align}\label{eq:epsilonZ}
\epsilon_{\mathrm{J}}=-(J^{\prime} \hat{H} J)^{-1} (J^{\prime} \hat{H} E \epsilon_{\mathrm{E}}+J^{\prime} c).
\end{align}
Inserting \eqref{eq:epsilonY} and \eqref{eq:epsilonZ} in \eqref{eq:epsilon} yields
\begin{align}
\begin{split}\label{eq:epsilonSolution}
\epsilon=&(E (G^{\mathcal{A}}E)^{-1} S^{\mathcal{A}}- J (J^{\prime}\hat{H}J)^{-1} J^{\prime} \hat{H} E (G^{\mathcal{A}}E)^{-1} S^{\mathcal{A}})x \\
&+(E-J (J^{\prime}\hat{H}J)^{-1} J^{\prime} \hat{H} E) (G^{\mathcal{A}}E)^{-1} W^{\mathcal{A}} \\ &- J (J^{\prime}\hat{H}J)^{-1} J^{\prime} c,
\end{split}
\end{align}
where $K_\epsilon$ and $b_\epsilon$ can be identified as in \eqref{eq:Kb}. Next, we multiply \eqref{eq:KKT_a_minmax} by $E^{\prime}$, i.e.
\begin{align}\label{eq:EqForLagrangeMinMax}
E^{\prime} \hat{H} \epsilon+E^{\prime} c+E^{\prime} (G^{\mathcal{A}})^{\prime} \lambda^{\mathcal{A}}  = 0.
\end{align}
Inserting \eqref{eq:epsilonSolution} in \eqref{eq:EqForLagrangeMinMax} and solving for $\lambda^{\mathcal{A}}$ results in 
\begin{align}\label{eq:lambdaSolution}
\lambda^{\mathcal{A}}=-\Theta^{\prime} E^{\prime} (\hat{H} K_{\epsilon}) x -\Theta^{\prime} E^{\prime} (\hat{H} b_{\epsilon}+c),
\end{align}
where $K_\lambda$ and $b_\lambda$ can be identified as in \eqref{eq:Kb}. The first part of the polytope $\mathcal{P}$, i.e., the first rows of $T$ and $d$ in \eqref{eq:Td}, can be computed by substituting the solution for $\epsilon$ acc. to \eqref{eq:epsilonSolution} into the inequality \eqref{eq:KKT_c_minmax}. The second part can be obtained by inserting the solution for $\lambda^{\mathcal{A}}$ acc. to \eqref{eq:lambdaSolution} in \eqref{eq:KKT_d_minmax}. \hfill $\square$
\end{proof}

By simple arithmetic operations the affine optimizer $\epsilon^\star(x)$ from Lemma 1 results in a robust optimal affine feedback law that yields the optimal control input $u^\star$ to system \eqref{eq:Sys} for all states $x \in \mathcal{P}$. This statement is based on the results in \cite{Ramirez2006} and \cite{Munoz2007} and is stated more precisely in Corollary \ref{cor:affineControlLaw}. In preparation we introduce
$K=K_{\mathcal{V}}^{\mathcal{M}}-K_{\infty}$, $K_{\mathcal{V}}=K_{\epsilon}^{\mathcal{N}} -H^{-1} F'$, $b=b_{\mathcal{V}}^{\mathcal{M}}$ and $b_{\mathcal{V}}=b_{\epsilon}^{\mathcal{N}}$ with $\mathcal{M}=\{1, \ldots, m\}$ and $\mathcal{N}=\{1, \ldots, mN\}$. 

\begin{corollary}\label{cor:affineControlLaw}
Let $x \in \mathcal{X}_f$ be arbitrary and $\mathcal{A}(x)$ the corresponding active set. Let the assumptions from Lemma \ref{lem:robustRegionalFeedbackLaw} hold. Then the robust optimal affine feedback law $u^\star(x)=K x + b$ yields the optimal control input $u^\star$ to system \eqref{eq:Sys} for all states $x \in \mathcal{P}$.
\end{corollary}

\begin{proof}
By combining \eqref{eq:VandZ}, \eqref{eq:augmentedVariable} and \eqref{eq:Kb} the optimal correction control input sequence can be expressed as 
\begin{align}\label{eq:affineV}
V^\star=K_{\mathcal{V}} x + b_{\mathcal{V}} 
\end{align}
with $K_{\mathcal{V}}=K_{\epsilon}^{\mathcal{N}} -H^{-1} F'$ and $b_{\mathcal{V}}=b_{\epsilon}^{\mathcal{N}}$. Substituting the first $m$ elements of \eqref{eq:affineV} in \eqref{eq:LQR_correction} results in the robust optimal affine feedback law 
\begin{align}
u^\star=K x + b
\end{align}
with $K=K_{\mathcal{V}}^{\mathcal{M}}-K_{\infty} $ and $b=b_{\mathcal{V}} ^{\mathcal{M}}$. \hfill $\square$
\end{proof}
\vspace{3mm}

Lemma \ref{lem:robustRegionalFeedbackLaw} and Corollary \ref{cor:affineControlLaw} can be used in a regional MPC approach as follows: As in \cite{Jost2015a}, Problem \eqref{eq:rereformulatedMPCProblem} is solved for the current state $x$, and the set of active constraints $\mathcal{A}(x)$ is determined according to \eqref{eq:sets}. Then a robust feedback law $u^\star(x)=K x+b$ and its polytope $\mathcal{P}$ are computed according to Lemma \ref{lem:robustRegionalFeedbackLaw} and Corollary \ref{cor:affineControlLaw}. The feedback law can be reused until its polytope has been left. When leaving the polytope, \eqref{eq:rereformulatedMPCProblem} must be solved again. If the rank condition on $G^{{\mathcal{A}}}$ for a state $x$ is not met, problem \eqref{eq:rereformulatedMPCProblem} must be solved in the next time step. We claim without giving details a rank violation rarely occurs.

\subsection{Extending min-max regional MPC by active set updates}\label{sec:AccActiveSetUp}
The new robust regional approach presented in the previous section can be extended by active set updates as proposed in \cite{Koenig2017a}. Note that the optimizer $\epsilon^\star(x)$ and the Lagrange multipliers $\lambda^\star(x)$ are piecewise affine functions of the parameter $x$ and $\epsilon^\star(x)$ is continuous and uniquely defined due to $H \succ 0$ (see \cite{Ramirez2006, Munoz2007, Tondel2003}). Thus, all requirements are met to apply the active set update approach (see \cite[Lemma 1]{Koenig2017a}). 
The extended approach can be described as follows: After leaving the current polytope the active set is updated along a line connecting the current and the previous state. This is done by analyzing the crossed facets of neighboring polytopes along the line. In doing so, several feedback laws and their polytopes along the line can be computed from the updated active sets without solving \eqref{eq:MPCProblem}. A new optimization problem has to be solved only if the active rows of the matrix $G$ on a crossed facet do not have full row rank or the weakly active set on a crossed facet consists of more than one constraint. For  more details we refer to \cite{Koenig2017a}.    

\begin{figure*}[t]
\centering
\captionsetup[subfloat]{justification=centering}
{
%
%
\begin{psfrags}%
\psfragscanon%
\scriptsize%
%
\psfrag{s01}[b][b]{\color[rgb]{0.15,0.15,0.15}\setlength{\tabcolsep}{0pt}\begin{tabular}{c}$x_2$\end{tabular}}%
\psfrag{s04}[t][t]{\color[rgb]{0.15,0.15,0.15}\setlength{\tabcolsep}{0pt}\begin{tabular}{c}$x_1$\end{tabular}}%
%
\color[rgb]{0.15,0.15,0.15}%
%
\psfrag{x01}[t][t]{-10}%
\psfrag{x02}[t][t]{-8}%
\psfrag{x03}[t][t]{-6}%
\psfrag{x04}[t][t]{-4}%
\psfrag{x05}[t][t]{-2}%
\psfrag{x06}[t][t]{0}%
\psfrag{x07}[t][t]{2}%
\psfrag{x08}[t][t]{4}%
\psfrag{x09}[t][t]{6}%
\psfrag{x10}[t][t]{8}%
\psfrag{x11}[t][t]{10}%
%
\psfrag{v01}[r][r]{-2}%
\psfrag{v02}[r][r]{-1}%
\psfrag{v03}[r][r]{0}%
\psfrag{v04}[r][r]{1}%
\psfrag{v05}[r][r]{2}%
\psfrag{v06}[r][r]{3}%
\psfrag{v07}[r][r]{4}%
\psfrag{v08}[r][r]{5}%
%
\includegraphics[width = 0.3\textwidth]{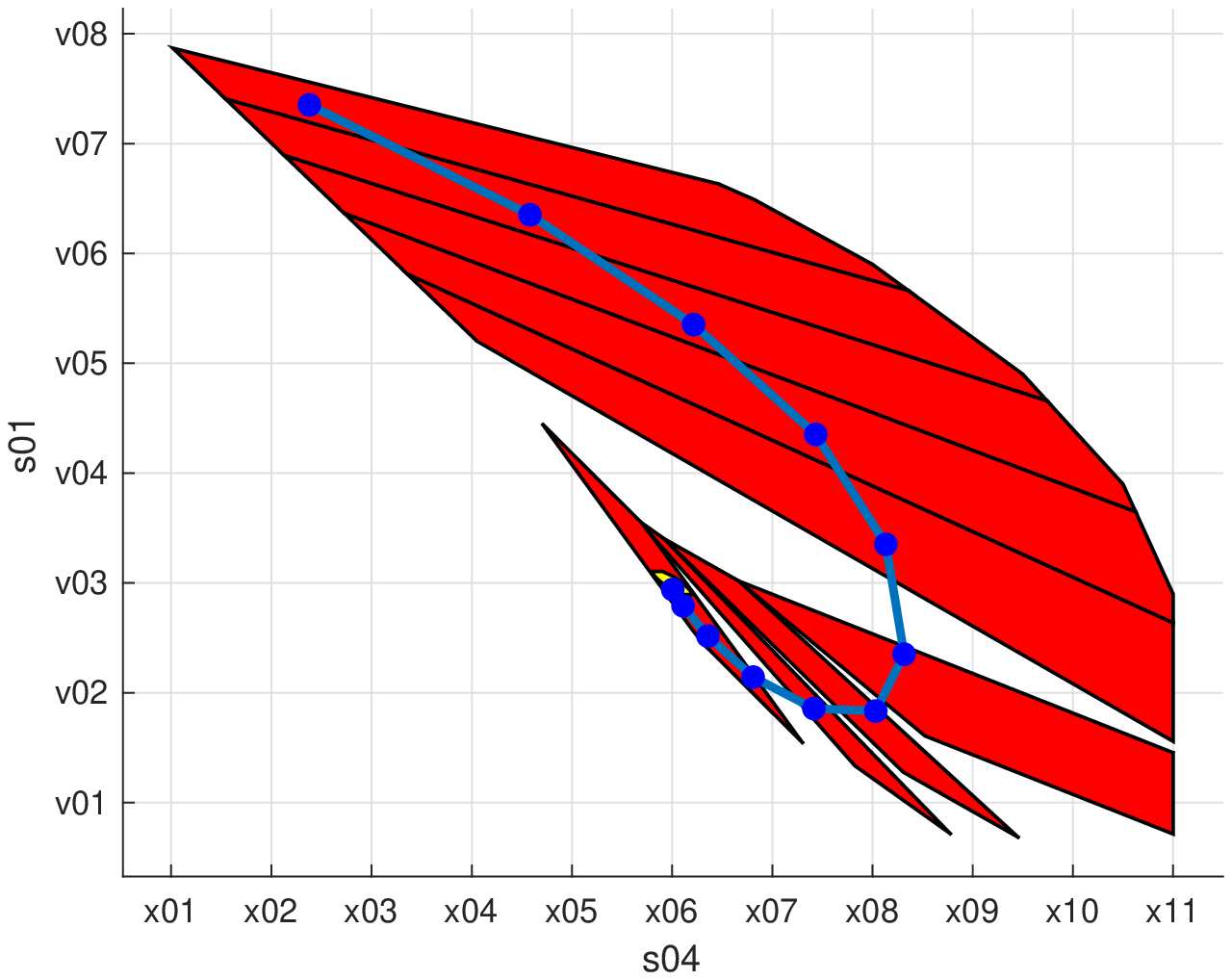}%
\end{psfrags}%
%
} \quad 
{
%
%
\begin{psfrags}%
\psfragscanon%
\scriptsize%
%
\psfrag{s02}[b][b]{\color[rgb]{0.15,0.15,0.15}\setlength{\tabcolsep}{0pt}\begin{tabular}{c}$x_2$\end{tabular}}%
\psfrag{s04}[t][t]{\color[rgb]{0.15,0.15,0.15}\setlength{\tabcolsep}{0pt}\begin{tabular}{c}$x_1$\end{tabular}}%
%
\color[rgb]{0.15,0.15,0.15}%
%
\psfrag{x01}[t][t]{-10}%
\psfrag{x02}[t][t]{-8}%
\psfrag{x03}[t][t]{-6}%
\psfrag{x04}[t][t]{-4}%
\psfrag{x05}[t][t]{-2}%
\psfrag{x06}[t][t]{0}%
\psfrag{x07}[t][t]{2}%
\psfrag{x08}[t][t]{4}%
\psfrag{x09}[t][t]{6}%
\psfrag{x10}[t][t]{8}%
\psfrag{x11}[t][t]{10}%
%
\psfrag{v01}[r][r]{-2}%
\psfrag{v02}[r][r]{-1}%
\psfrag{v03}[r][r]{0}%
\psfrag{v04}[r][r]{1}%
\psfrag{v05}[r][r]{2}%
\psfrag{v06}[r][r]{3}%
\psfrag{v07}[r][r]{4}%
\psfrag{v08}[r][r]{5}%
%
\includegraphics[width = 0.3\textwidth]{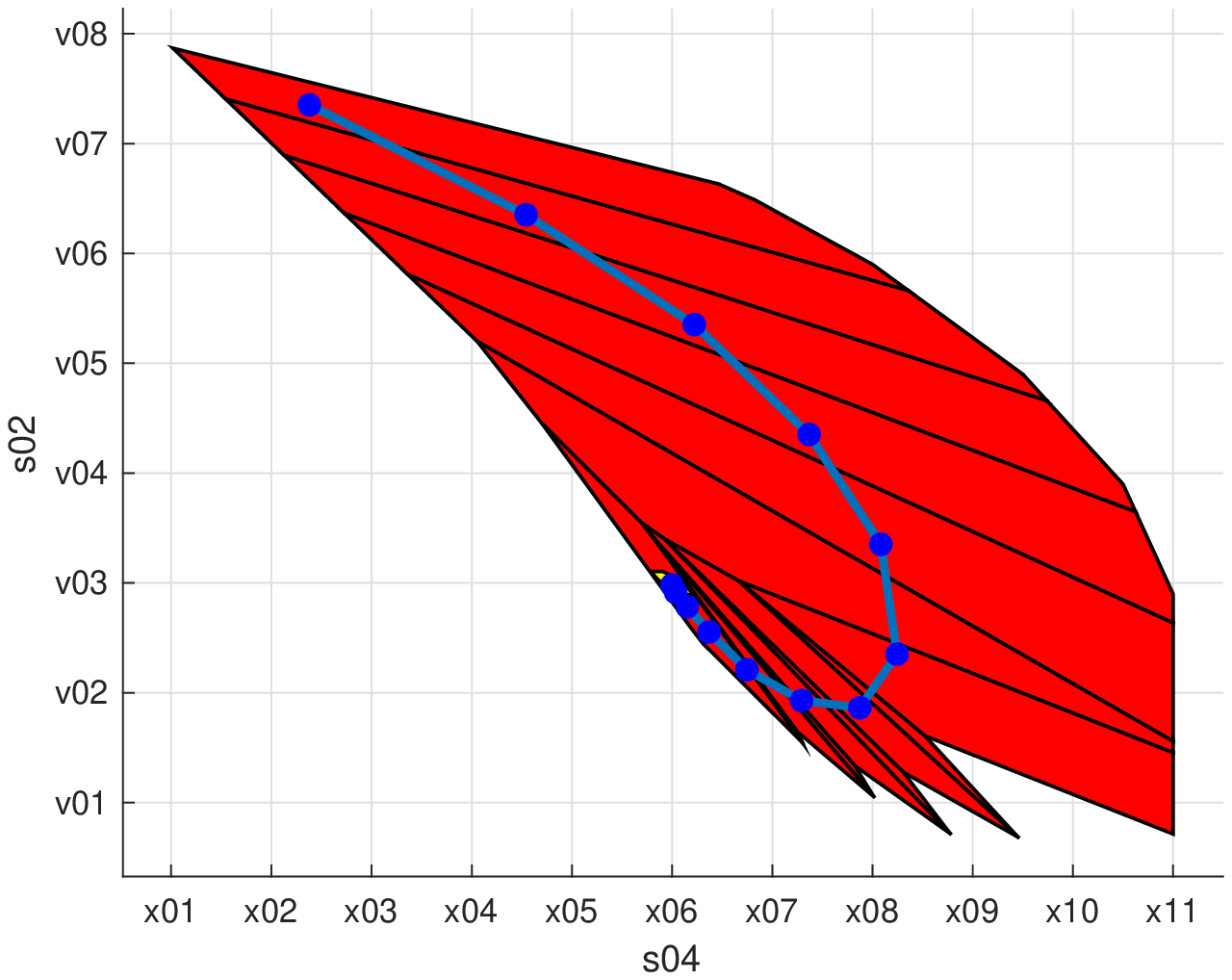}%
\end{psfrags}%
%
} \quad 
{
%
%
\begin{psfrags}%
\psfragscanon%
\scriptsize%
%
\psfrag{s01}[t][t]{\color[rgb]{0.15,0.15,0.15}\setlength{\tabcolsep}{0pt}\begin{tabular}{c}$x_1$\end{tabular}}%
\psfrag{s03}[b][b]{\color[rgb]{0.15,0.15,0.15}\setlength{\tabcolsep}{0pt}\begin{tabular}{c}$x_2$\end{tabular}}%
%
\color[rgb]{0.15,0.15,0.15}%
%
\psfrag{x01}[t][t]{-10}%
\psfrag{x02}[t][t]{-8}%
\psfrag{x03}[t][t]{-6}%
\psfrag{x04}[t][t]{-4}%
\psfrag{x05}[t][t]{-2}%
\psfrag{x06}[t][t]{0}%
\psfrag{x07}[t][t]{2}%
\psfrag{x08}[t][t]{4}%
\psfrag{x09}[t][t]{6}%
\psfrag{x10}[t][t]{8}%
\psfrag{x11}[t][t]{10}%
%
\psfrag{v01}[r][r]{-4}%
\psfrag{v02}[r][r]{-3}%
\psfrag{v03}[r][r]{-2}%
\psfrag{v04}[r][r]{-1}%
\psfrag{v05}[r][r]{0}%
\psfrag{v06}[r][r]{1}%
\psfrag{v07}[r][r]{2}%
\psfrag{v08}[r][r]{3}%
\psfrag{v09}[r][r]{4}%
\psfrag{v10}[r][r]{5}%
%
\includegraphics[width = 0.3\textwidth]{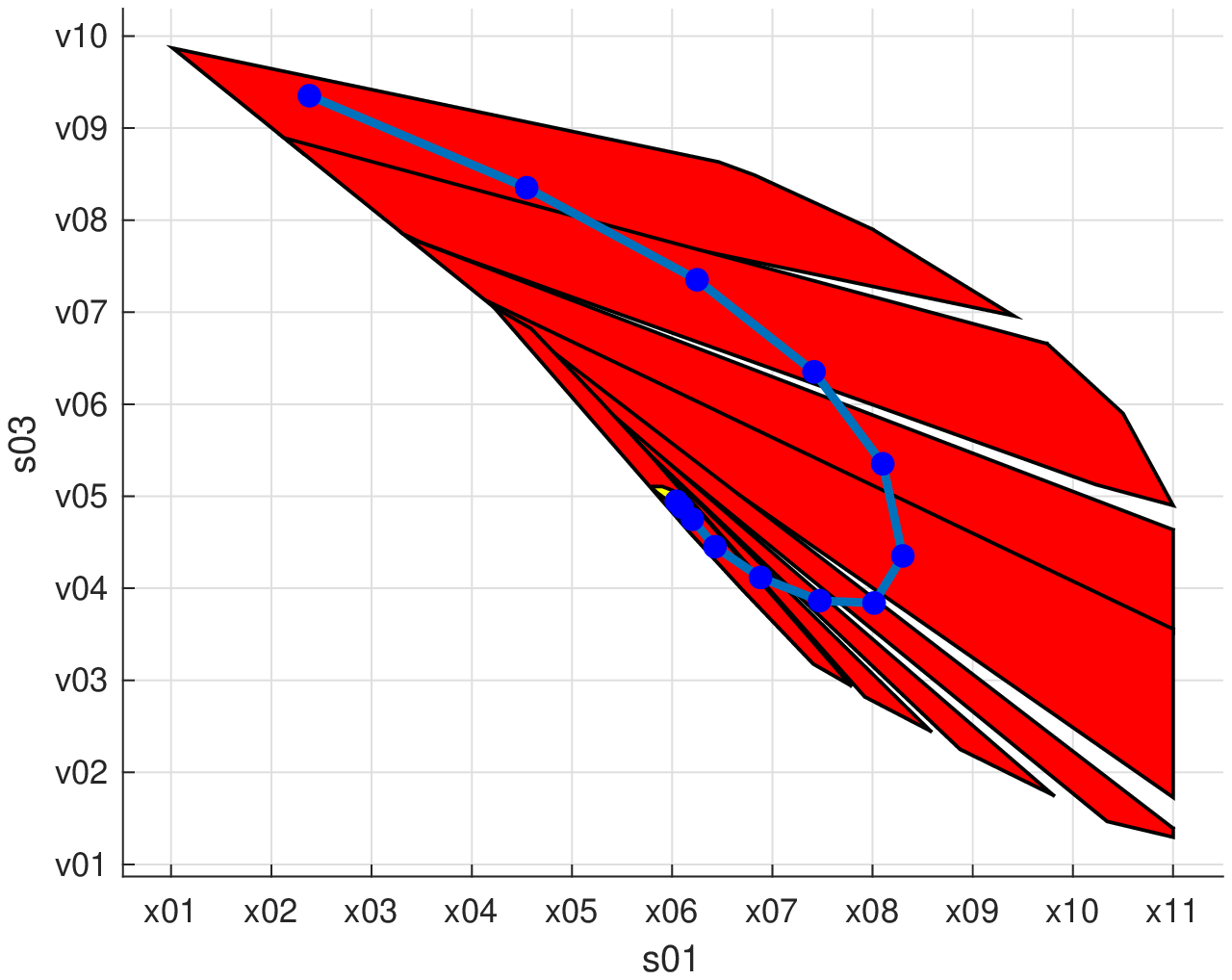}%
\end{psfrags}%
%
} \\ \vspace{5pt}
\subfloat[min-max regional MPC]{
%
%
\begin{psfrags}%
\psfragscanon%
\scriptsize%
%
\psfrag{s08}[b][b]{\color[rgb]{0.15,0.15,0.15}\setlength{\tabcolsep}{0pt}\begin{tabular}{c}$u(k)$\end{tabular}}%
\psfrag{s09}[b][b]{\color[rgb]{0.15,0.15,0.15}\setlength{\tabcolsep}{0pt}\begin{tabular}{c}$w(k)$\end{tabular}}%
\psfrag{s12}[b][b]{\color[rgb]{0.15,0.15,0.15}\setlength{\tabcolsep}{0pt}\begin{tabular}{c}$e(k)$\end{tabular}}%
\psfrag{s13}[t][t]{\color[rgb]{0.15,0.15,0.15}\setlength{\tabcolsep}{0pt}\begin{tabular}{c}$k$\end{tabular}}%
\psfrag{s15}[b][b]{\color[rgb]{0.15,0.15,0.15}\setlength{\tabcolsep}{0pt}\begin{tabular}{c}$x(k)$\end{tabular}}%
%
\color[rgb]{0.15,0.15,0.15}%
%
\psfrag{x01}[t][t]{0}%
\psfrag{x02}[t][t]{1}%
\psfrag{x03}[t][t]{2}%
\psfrag{x04}[t][t]{3}%
\psfrag{x05}[t][t]{4}%
\psfrag{x06}[t][t]{5}%
\psfrag{x07}[t][t]{6}%
\psfrag{x08}[t][t]{7}%
\psfrag{x09}[t][t]{8}%
\psfrag{x10}[t][t]{9}%
\psfrag{x11}[t][t]{10}%
\psfrag{x12}[t][t]{0}%
\psfrag{x13}[t][t]{1}%
\psfrag{x14}[t][t]{2}%
\psfrag{x15}[t][t]{3}%
\psfrag{x16}[t][t]{4}%
\psfrag{x17}[t][t]{5}%
\psfrag{x18}[t][t]{6}%
\psfrag{x19}[t][t]{7}%
\psfrag{x20}[t][t]{8}%
\psfrag{x21}[t][t]{9}%
\psfrag{x22}[t][t]{10}%
\psfrag{x23}[t][t]{11}%
\psfrag{x24}[t][t]{0}%
\psfrag{x25}[t][t]{1}%
\psfrag{x26}[t][t]{2}%
\psfrag{x27}[t][t]{3}%
\psfrag{x28}[t][t]{4}%
\psfrag{x29}[t][t]{5}%
\psfrag{x30}[t][t]{6}%
\psfrag{x31}[t][t]{7}%
\psfrag{x32}[t][t]{8}%
\psfrag{x33}[t][t]{9}%
\psfrag{x34}[t][t]{10}%
\psfrag{x35}[t][t]{11}%
\psfrag{x36}[t][t]{0}%
\psfrag{x37}[t][t]{1}%
\psfrag{x38}[t][t]{2}%
\psfrag{x39}[t][t]{3}%
\psfrag{x40}[t][t]{4}%
\psfrag{x41}[t][t]{5}%
\psfrag{x42}[t][t]{6}%
\psfrag{x43}[t][t]{7}%
\psfrag{x44}[t][t]{8}%
\psfrag{x45}[t][t]{9}%
\psfrag{x46}[t][t]{10}%
\psfrag{x47}[t][t]{11}%
%
\psfrag{v01}[r][r]{-1}%
\psfrag{v02}[r][r]{0}%
\psfrag{v03}[r][r]{1}%
\psfrag{v04}[r][r]{0}%
\psfrag{v05}[r][r]{}%
\psfrag{v06}[r][r]{1}%
\psfrag{v07}[r][r]{-2}%
\psfrag{v08}[r][r]{0}%
\psfrag{v09}[r][r]{2}%
\psfrag{v10}[r][r]{-10}%
\psfrag{v11}[r][r]{0}%
\psfrag{v12}[r][r]{10}%
%
\includegraphics[width = 0.3\textwidth, height = 0.25\textwidth]{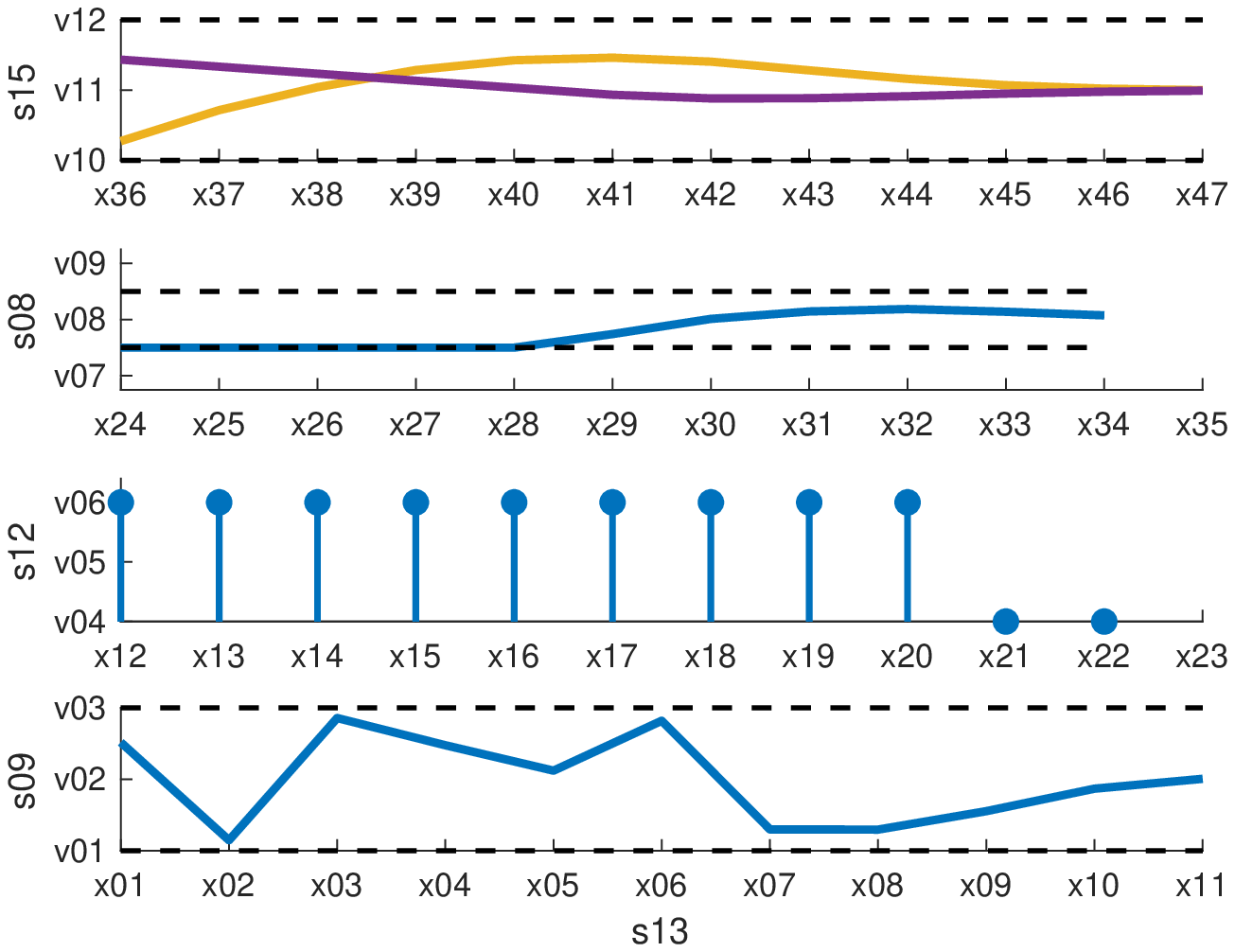}%
\end{psfrags}%
%
}\quad \quad
\subfloat[using active set updates]{
%
%
\begin{psfrags}%
\psfragscanon%
\scriptsize%
%
\psfrag{s05}[b][b]{\color[rgb]{0.15,0.15,0.15}\setlength{\tabcolsep}{0pt}\begin{tabular}{c}$x(k)$\end{tabular}}%
\psfrag{s07}[b][b]{\color[rgb]{0.15,0.15,0.15}\setlength{\tabcolsep}{0pt}\begin{tabular}{c}$e(k)$\end{tabular}}%
\psfrag{s09}[t][t]{\color[rgb]{0.15,0.15,0.15}\setlength{\tabcolsep}{0pt}\begin{tabular}{c}$k$\end{tabular}}%
\psfrag{s10}[b][b]{\color[rgb]{0.15,0.15,0.15}\setlength{\tabcolsep}{0pt}\begin{tabular}{c}$w(k)$\end{tabular}}%
\psfrag{s12}[b][b]{\color[rgb]{0.15,0.15,0.15}\setlength{\tabcolsep}{0pt}\begin{tabular}{c}$u(k)$\end{tabular}}%
%
\color[rgb]{0.15,0.15,0.15}%
%
\psfrag{x01}[t][t]{0}%
\psfrag{x02}[t][t]{1}%
\psfrag{x03}[t][t]{2}%
\psfrag{x04}[t][t]{3}%
\psfrag{x05}[t][t]{4}%
\psfrag{x06}[t][t]{5}%
\psfrag{x07}[t][t]{6}%
\psfrag{x08}[t][t]{7}%
\psfrag{x09}[t][t]{8}%
\psfrag{x10}[t][t]{9}%
\psfrag{x11}[t][t]{10}%
\psfrag{x12}[t][t]{11}%
\psfrag{x13}[t][t]{0}%
\psfrag{x14}[t][t]{2}%
\psfrag{x15}[t][t]{4}%
\psfrag{x16}[t][t]{6}%
\psfrag{x17}[t][t]{8}%
\psfrag{x18}[t][t]{10}%
\psfrag{x19}[t][t]{12}%
\psfrag{x20}[t][t]{0}%
\psfrag{x21}[t][t]{2}%
\psfrag{x22}[t][t]{4}%
\psfrag{x23}[t][t]{6}%
\psfrag{x24}[t][t]{8}%
\psfrag{x25}[t][t]{10}%
\psfrag{x26}[t][t]{12}%
\psfrag{x27}[t][t]{0}%
\psfrag{x28}[t][t]{2}%
\psfrag{x29}[t][t]{4}%
\psfrag{x30}[t][t]{6}%
\psfrag{x31}[t][t]{8}%
\psfrag{x32}[t][t]{10}%
\psfrag{x33}[t][t]{12}%
%
\psfrag{v01}[r][r]{-1}%
\psfrag{v02}[r][r]{0}%
\psfrag{v03}[r][r]{1}%
\psfrag{v04}[r][r]{0}%
\psfrag{v05}[r][r]{}%
\psfrag{v06}[r][r]{1}%
\psfrag{v07}[r][r]{-2}%
\psfrag{v08}[r][r]{0}%
\psfrag{v09}[r][r]{2}%
\psfrag{v10}[r][r]{-10}%
\psfrag{v11}[r][r]{0}%
\psfrag{v12}[r][r]{10}%
%
\includegraphics[width = 0.3\textwidth, height = 0.25\textwidth]{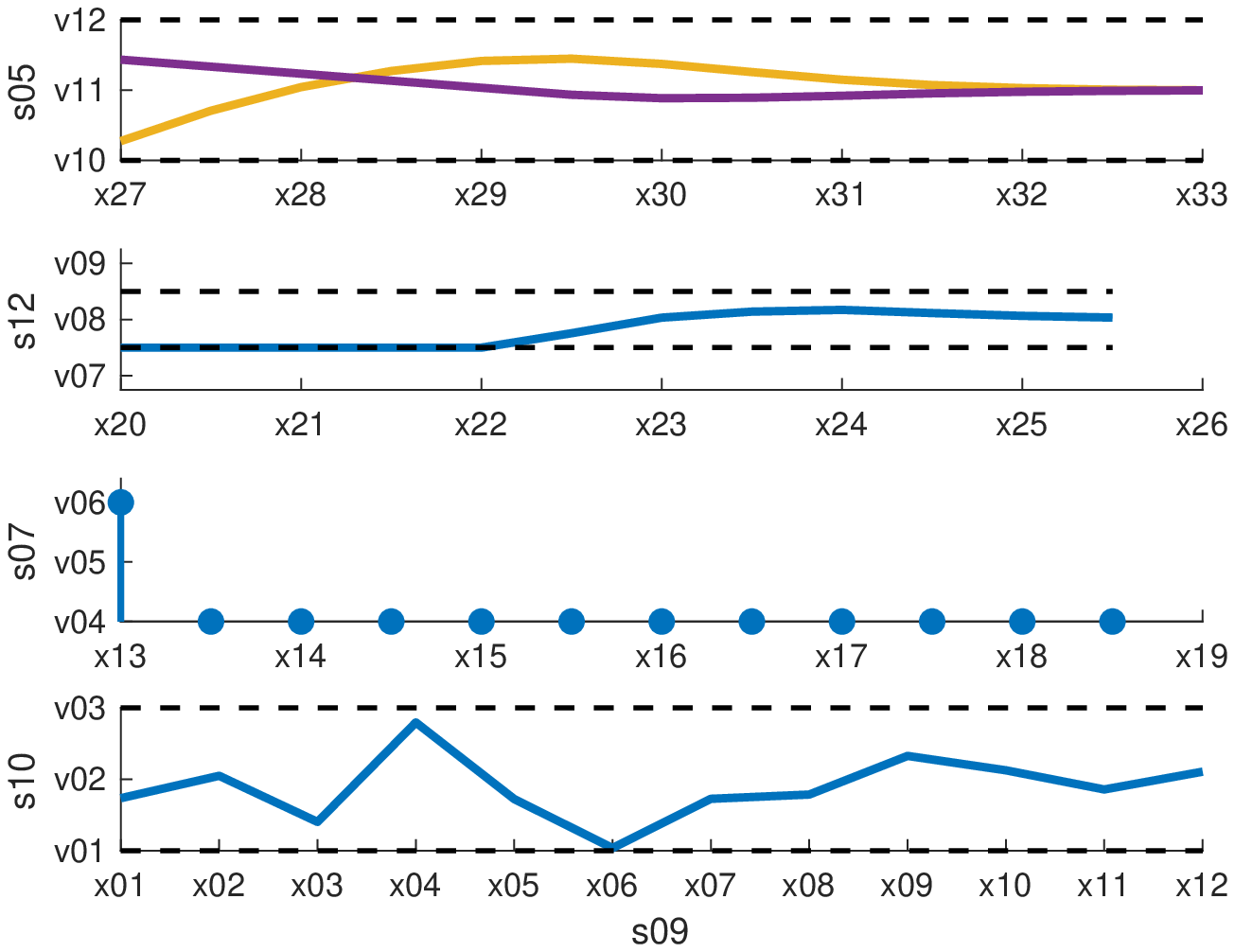}%
\end{psfrags}%
%
}\quad \quad
\subfloat[using suboptimality]{
%
%
\begin{psfrags}%
\psfragscanon%
\scriptsize%
%
\psfrag{s02}[t][t]{\color[rgb]{0.15,0.15,0.15}\setlength{\tabcolsep}{0pt}\begin{tabular}{c}$k$\end{tabular}}%
\psfrag{s04}[b][b]{\color[rgb]{0.15,0.15,0.15}\setlength{\tabcolsep}{0pt}\begin{tabular}{c}$u(k)$\end{tabular}}%
\psfrag{s06}[b][b]{\color[rgb]{0.15,0.15,0.15}\setlength{\tabcolsep}{0pt}\begin{tabular}{c}$x(k)$\end{tabular}}%
\psfrag{s09}[b][b]{\color[rgb]{0.15,0.15,0.15}\setlength{\tabcolsep}{0pt}\begin{tabular}{c}$w(k)$\end{tabular}}%
\psfrag{s16}[b][b]{\color[rgb]{0.15,0.15,0.15}\setlength{\tabcolsep}{0pt}\begin{tabular}{c}$e(k)$\end{tabular}}%
%
\color[rgb]{0.15,0.15,0.15}%
%
\psfrag{x01}[t][t]{0}%
\psfrag{x02}[t][t]{1}%
\psfrag{x03}[t][t]{2}%
\psfrag{x04}[t][t]{3}%
\psfrag{x05}[t][t]{4}%
\psfrag{x06}[t][t]{5}%
\psfrag{x07}[t][t]{6}%
\psfrag{x08}[t][t]{7}%
\psfrag{x09}[t][t]{8}%
\psfrag{x10}[t][t]{9}%
\psfrag{x11}[t][t]{10}%
\psfrag{x12}[t][t]{11}%
\psfrag{x13}[t][t]{0}%
\psfrag{x14}[t][t]{2}%
\psfrag{x15}[t][t]{4}%
\psfrag{x16}[t][t]{6}%
\psfrag{x17}[t][t]{8}%
\psfrag{x18}[t][t]{10}%
\psfrag{x19}[t][t]{12}%
\psfrag{x20}[t][t]{0}%
\psfrag{x21}[t][t]{2}%
\psfrag{x22}[t][t]{4}%
\psfrag{x23}[t][t]{6}%
\psfrag{x24}[t][t]{8}%
\psfrag{x25}[t][t]{10}%
\psfrag{x26}[t][t]{12}%
\psfrag{x27}[t][t]{0}%
\psfrag{x28}[t][t]{2}%
\psfrag{x29}[t][t]{4}%
\psfrag{x30}[t][t]{6}%
\psfrag{x31}[t][t]{8}%
\psfrag{x32}[t][t]{10}%
\psfrag{x33}[t][t]{12}%
%
\psfrag{v01}[r][r]{-1}%
\psfrag{v02}[r][r]{0}%
\psfrag{v03}[r][r]{1}%
\psfrag{v04}[r][r]{0}%
\psfrag{v05}[r][r]{}%
\psfrag{v06}[r][r]{1}%
\psfrag{v07}[r][r]{-2}%
\psfrag{v08}[r][r]{0}%
\psfrag{v09}[r][r]{2}%
\psfrag{v10}[r][r]{-10}%
\psfrag{v11}[r][r]{0}%
\psfrag{v12}[r][r]{10}%
%
\includegraphics[width = 0.3\textwidth, height = 0.25\textwidth]{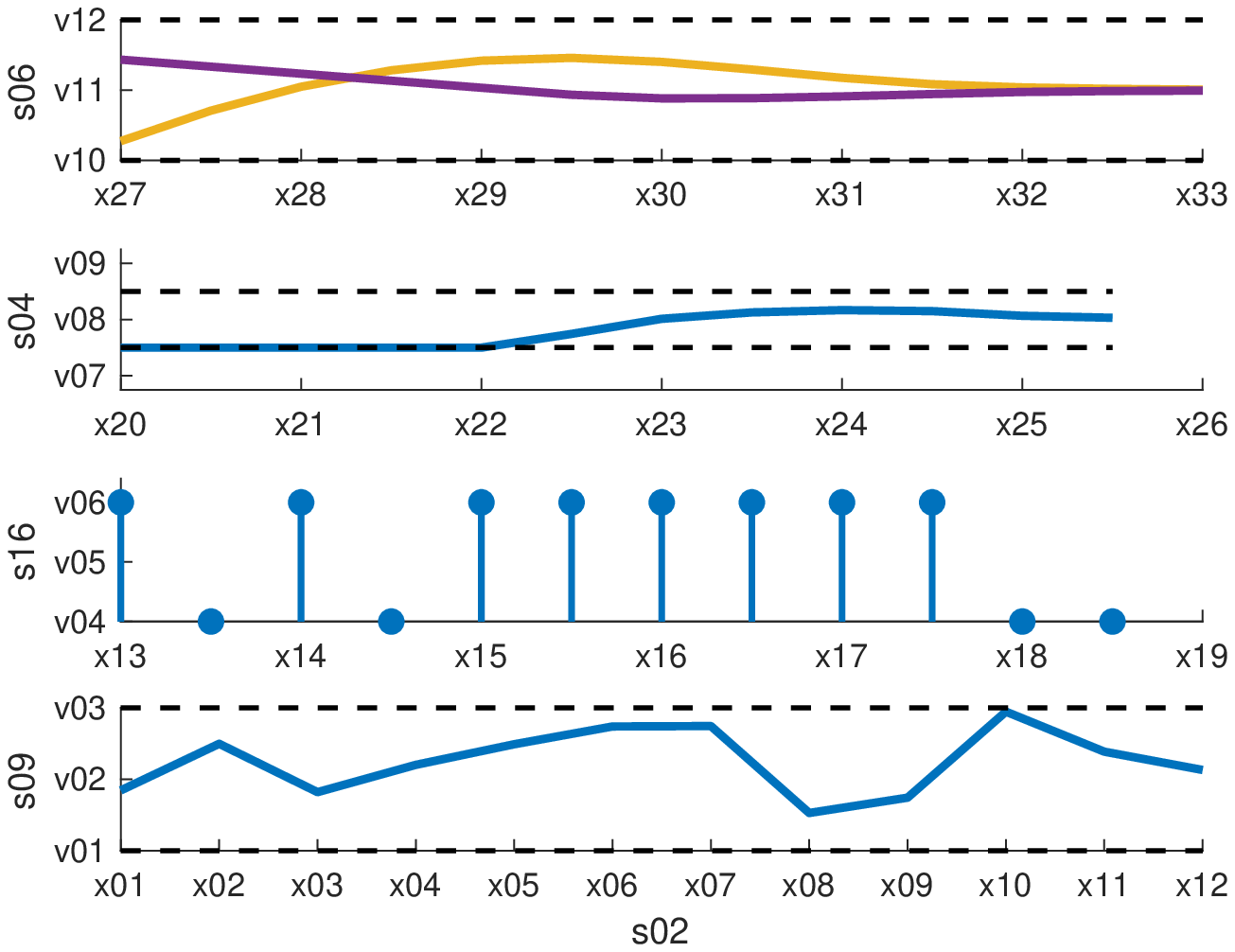}%
\end{psfrags}%
%
}\\
\caption{The figure shows closed-loop trajectories and polytopes (red) in state-space (upper plots) and the time series of the states $x(k)$, inputs $u(k)$, disturbances $w(k)$ and the event $e(k)$ (bottom plots) for the new min-max regional MPC approach (a) and the extensions using active set updates (b) and suboptimality (c). } 
\label{fig:BeispielSISO}
\end{figure*}  

\setlength{\tabcolsep}{2pt}
\begin{table*}[t]
\footnotesize
{\normalsize Table 1. The table compares the new min-max, the tube-based and the classical regional MPC approach. It also shows the results for the approaches extended by active set updates and suboptimality. The suboptimal approaches have been executed with $M=15$ steps (see Section \ref{sec:AccSuboptimality} for details).}
\begin{center}
\begin{tabular}{|c|c|c|c|c|c|c|c|c|c|c|c|c|c|c|c|c|}
\hline
 \multirow{2}{*}{approach} & \multirow{2}{*}{ext.} & \multicolumn{5}{|c|}{$N=3$}    & \multicolumn{5}{|c|}{$N=5$}  &  \multicolumn{5}{|c|}{$N=10$} \\
\cline{3-17}
     & & $q$ & $p$ & steps & $\Delta$QP (in \%) & $\Delta T$ & $q$ & $p$ & steps & $\Delta$QP (in \%) & $\Delta T$ &$q$&$p$& steps& $\Delta$QP (in \%)& $\Delta T$   \\ \hline \hline
      \multirow{3}{*}{\parbox{1.9cm}{\centering min-max regional MPC}} & basic & \multirow{3}{*}{32} & \multirow{3}{*}{4} & 8.65& -3.15 (-36.36) &-23.18& \multirow{3}{*}{68} & \multirow{3}{*}{6} &9.35 & -2.10 (-22.49) & -13.18  & \multirow{3}{*}{1090} & \multirow{3}{*}{11} &10.17&-1.27 (-12.53)& -6.61\\ 
   & active set up.&&&8.35&-7.31 (-87.48)&-58.18&&& 9.36 & -8.32 (-88.85) & -57.23 &&&9.92&-8.86 (-89.32)&-83.85\\ 
     & suboptimality &&&8.01&-3.33 (-41.60)&-35.00&&& 9.23 & -3.08 (-33.41) & -25.72 &&&9.96&-2.16 (-21.72)&-15.16\\ \hline 
     \multirow{3}{*}{\parbox{1.9cm}{\centering tube-based regional MPC}} & basic &\multirow{3}{*}{74}&\multirow{3}{*}{5}& 9.76 &-0.80 (-8.18)& +1.57 & \multirow{3}{*}{86}&\multirow{3}{*}{7}& 10.55 & 0.82 (-7.73) & +0.18 & \multirow{3}{*}{116} & \multirow{3}{*}{12} & 11.19&-0.89 (-7.99)& -0.38 \\ 
   &  active set up. & &&9.76&-8.75 (-89.67)&-33.18&&&10.46 & -9.46 (-90.39) & -40.60 &&&11.15&-10.14 (-90.99)&-49.11\\ 
     & suboptimality &&&9.50&-7.21 (-75.92)&-70.85&&& 10.38 & -7.71 (-74.24) & -68.54 &&&11.01& -7.93 (-72.09)& -67.22\\ \hline
\multirow{3}{*}{\parbox{1.9cm}{\centering regional MPC (undisturbed)}} & basic &\multirow{3}{*}{24}&\multirow{3}{*}{3}& 2.00 & 0.00 (-0.00) & +12.70 &\multirow{3}{*}{36}&\multirow{3}{*}{5}& 2.78 & -0.00 (-0.00) & +15.56 & \multirow{3}{*}{66} &\multirow{3}{*}{10} &3.11& -0.00 (-0.00)& +14.85  \\ 
   & active set up. &&&2.00&1.00 (-50.05)&-34.92&&& 2.78 & -1.78 (-64.08)& -47.78 & &&3.11& -2.11 (-67.87)& -48.51\\ 
     & suboptimality& &&2.00&0.72 (-36.16)&-23.81&&& 2.78 & -1.10 (-39.37) & -26.67& &&3.11& -1.24 (-39.78)& -28.71\\ \hline 
\end{tabular}
\label{tab:Results}
\end{center}
\end{table*}

\subsection{Extending min-max regional MPC by suboptimality}\label{sec:AccSuboptimality}
The new robust regional MPC approach can also be extended by suboptimality as proposed in \cite{SchulzeDarup2017CDC}. The region $\mathcal{P}$ in Lemma \ref{lem:robustRegionalFeedbackLaw} and Corollary \ref{cor:affineControlLaw} is the intersection of two regions $\mathcal{C}$ (defined by the first rows in \eqref{eq:Td}) and $\mathcal{O}$ (defined by the last rows in \eqref{eq:Td}) that guarantee constraint satisfaction and optimality of the input resulting from $K x + b$. If $x \notin \mathcal{P}=\mathcal{C} \cap \mathcal{O}$ but $x \in \mathcal{C}$ then $Kx + b$ yields not an optimal input but a feasible one, i.e., an input that satisfies the constraints in \eqref{eq:constraints}. 
The idea is to reuse a feedback law even outside its polytope $\mathcal{P}$ of optimality as long as it is feasible and stabilizes the system. For that, a feedback law is reused as long as the system is located inside $\mathcal{C}$. Stability can be guaranteed by switching irreversibly to a point-by-point solution of problem \eqref{eq:MPCProblem} if the RPI set is not reached after a user-defined number of $M$ steps. Note that allowing suboptimality enlarges the validity region for a feedback law resulting in a reduced number of online optimizations and thus a reduced overall computational effort. For more details, we refer to \cite{SchulzeDarup2017CDC}.

\section{Numerical example}\label{sec:Example}
We investigate the computational effort of the approach proposed in the previous section and compare it to both an existing tube-based approach, and the original undisturbed approach from \cite{Jost2015a}. We consider a double integrator taken from \cite{Munoz2007} that results in a system of the form \eqref{eq:Sys} with
	\begin{align}
		A=\begin{pmatrix}  
		1 & 1 \\
		0 & 1
		\end{pmatrix}, \quad
		B=\begin{pmatrix}
		0\\
		1
		\end{pmatrix},\quad
		D=\begin{pmatrix}
		0.1\\
		0
		\end{pmatrix}. \nonumber
	\end{align}
	The system must respect the state and input constraints $-10 \leq x_i \leq 10$, $i=1,2$ and $-1 \leq u_1 \leq 1$, respectively, and the uncertainty is bounded by $-1 \leq w_1 \leq 1$. We choose the weighting matrices $Q=\text{diag}(1,1)$ and $R=10$. The weighting $P$ is chosen as the solution of the discrete-time algebraic Riccati equation. The RPI set $\mathcal{R}$ is computed with the procedure from \cite{Rakovic2005}. 
The terminal set $\mathcal{T}$ is computed with the procedure from \cite{Gilbert1991}. For the tube-based and min-max approach, we use a robust terminal set according to \cite{MAYNE2005}. The sets $\mathcal{R}$ and $\mathcal{T}$ are given in the appendix.  
We vary the horizon between $N=3$, $N=5$ and $N=10$ resulting in QPs of increasing size. The number of constraints $q$ and optimization variables $p$ of the QPs solved in the different approaches are given in Table 1.

\subsection{Analysis of a sample trajectory}
Figure \ref{fig:BeispielSISO} shows the results of the controlled system for an initial state for $N=5$. The proposed min-max regional MPC approach (a) and its extensions using active set updates (b) and suboptimality (c) are compared to each other. The upper plots show the state-space trajectories along with the computed polytopes (red). The bottom plots show the trajectories of the states $x(k)$, the inputs $u(k)$, the disturbances $w(k)$ and the indicator function $e(k)$ with $e(k)=1$ if a QP is solved and $e(k)=0$ otherwise. Constraints are depicted with black dashed lines. 

All three approaches regulate the system into the robust positively invariant set around the origin as expected. The state-space trajectories are very similar despite varying disturbance sequences. 
With the basic approach (a) nine QPs and thus nine polytopes and affine feedback laws must be calculated. 
Using active set updates (b) results in 15 polytopes and feedback laws by solving only a single QP. The same polytopes as in (a) are calculated but also polytopes in between. 
By exploiting suboptimality (c) the system can be regulated to the RPI set with eight affine feedback laws and polytopes. Some of the polytopes are larger than the polytopes in (a) due to suboptimality. Because of this, a QP can be saved in time steps one and three.

\subsection{Statistical analysis}

Table 1 shows the results from a statistical analysis. 
For sake of comparison, we also give the results for the tube-based regional MPC approach from \cite{SchulzeDarup2017CDC} and the classical regional MPC approach from \cite{Jost2015a}, both extended by active set updates and suboptimality. We refer to the corresponding publications \cite{Jost2015a}, \cite{Koenig2017a} and \cite{SchulzeDarup2017CDC} for more details. For each approach we have chosen 500 random initial states from the corresponding set $\mathcal{X}_f \setminus \mathcal{T}$. 

For every control scheme, we simulated trajectories of the controlled system until we reached $\mathcal{T}$ (regional MPC approach) respectively $\mathcal{R}$ (tube-based and min-max regional MPC approach). We used the same set $\mathcal{R}$ for both robust approaches. We applied randomly chosen disturbances to the system while simulating the tube-based and the min-max approach. The classical regional MPC approach was simulated without disturbances. The table states the average number of steps up to the set $\mathcal{T}$ (respectively $\mathcal{R}$) and the average number of avoided QPs ($\Delta$QP) per trajectory.  Moreover, the columns $\Delta T$ state the reduction of the overall computational effort (in percent) related to a classical point-by-point solution of the underlying optimization problem. The effort was determined by measuring the \textsc{Matlab} execution time. Note that $\Delta T$ takes not only the effort for solving the QPs but also the additional effort for computing the regional feedback laws into account.  
We stopped the simulation in $\mathcal{T}$ (respectively $\mathcal{R}$) since the unconstrained LQR can be used inside these sets without solving additional QPs.  

We analyse the results exemplary for the case $N=5$ 
in more detail (see Table 1). Note that the results for the cases $N=3$ and $N=10$ are similar to the case $N=5$ and are summarized at the end of this section. 
\begin{itemize}
\item In the new min-max regional MPC approach (rows 3-5) the system requires $9.3$ steps on average to reach the set $\mathcal{R}$. With the basic approach (see Fig. \ref{fig:BeispielSISO} (a)) $\unit[22.49]{\%}$ of the QPs can be avoided along the trajectory. This results in a decrease of the computational effort by $\unit[13.18]{\%}$. Exploiting suboptimality (see Fig. \ref{fig:BeispielSISO} (c)) decreases the number of QPs and the computational effort by $\unit[33.41]{\%}$ and $\unit[25.72]{\%}$, respectively. With active set updates (see Fig. \ref{fig:BeispielSISO} (b)) $\unit[88.85]{\%}$ of the QPs are avoided and the computational effort decreases by $\unit[57.23]{\%}$. 
\item The tube-based regional MPC approaches (rows 6-8) need $10.5$ steps on average to reach the set $\mathcal{R}$. With the basic approach $\unit[7.73]{\%}$ of the QPs are avoided. The computational effort is slightly increased by $\unit[0.18]{\%}$, however, because of the additional effort for computing the robust feedback laws. By exploiting suboptimality in the tube-based approach the number of QPs is decreased by $\unit[74.24]{\%}$ and the computational effort is reduced by $\unit[68.54]{\%}$. Tube-based MPC with active set updates avoids $\unit[90.39]{\%}$ of the QPs and the computational effort is reduced by $\unit[40.60]{\%}$.
\item In the classical regional MPC approach (without disturbances, rows 9-11) the number of QPs cannot be reduced with the basic approach.
Moreover, the computational effort slightly increases because of the additional computations for the feedback laws. Exploiting suboptimality leads to a decrease in the number of QPs and the computational effort by $\unit[39.37]{\%}$ and $\unit[26.67]{\%}$, respectively. Using active set updates leads to a decrease in the number of QPs and the computational effort by $\unit[64.08]{\%}$ and $\unit[47.78]{\%}$, respectively. 
\end{itemize}

The results for the horizons $N=3$ and $N=10$ are similar to the case $N=5$ and can be summarized as follows: With the regional MPC idea from \cite{Jost2015a} the computational effort for solving min-max MPC problems can be reduced. The highest savings can be achieved by using active set updates, which lead to a reduction of $\unit[83.85]{\%}$ for a problem with $N=10$ and $1090$ constraints. Moreover, the basic min-max regional MPC approach has greater savings than the basic tube-based regional MPC approach. However, exploiting suboptimality yields better results in the tube-based approach. 

\subsection{Choosing a robust approach based on the horizon}\label{sec:NumberConstraints}
We compute the horizon up to which the application of the min-max approach is still appropriate. 
The number of constraints in the min-max approach is 
\begin{align}
(2s)^N+2N(m+n)+q_{\mathcal{T}_{M}}, \nonumber
\end{align}
where the first term belongs to the constraints associated with the vertices $\text{Ver}(\mathcal{D}_{\mathrm{N}})$ with the disturbance dimension $s$, the second one belongs to state and input constraints and the last one belongs to the $q_{\mathcal{T}_{M}}$ terminal constraints.
The number of constraints in the tube-based approach (see \cite{SchulzeDarup2017CDC} for details) is 
\begin{align}
q_{\mathcal{R}}+2N(m+n)+q_{\mathcal{T}_{T}}, \nonumber
\end{align}
where the first term belongs to the $q_{\mathcal{R}}$ initial constraints (hyperplanes of $\mathcal{R}$), the second one belongs to state and input constraints and the last one belongs to the $q_{\mathcal{T}_{T}}$ terminal constraints.
The min-max approach has fewer constraints than the tube-based approach as long as
\begin{align} \label{eq:conditionConstraints}
(2s)^N+q_{\mathcal{T}_{M}} < q_{\mathcal{R}}+q_{\mathcal{T}_{T}}.     
\end{align}
Solving \eqref{eq:conditionConstraints} for $N$ results in
\begin{align}\label{eq:conditionHorizon}
N < \frac{\log(q_\mathcal{R}+q_{\mathcal{T}_{T}}-q_{\mathcal{T}_{M}})}{\log(2s)}.
\end{align}
With \eqref{eq:conditionHorizon} a horizon $\hat{N}$ can be computed up to which the min-max problem \eqref{eq:rereformulatedMPCProblem} has fewer constraints than the tube-based problem proposed in \cite{SchulzeDarup2017CDC}. 

For the considered double integrator $\hat{N}=5$ holds. Note that for horizon $N=10$ the min-max problem has almost ten times as many constraints as the tube-based problem (see $q$ in Table 1). Note that the optimization effort increases with the number of constraints and thus the tube-based approach should be preferred for $N=10$.

\section{Conclusions}\label{sec:Conclusion}
We applied the regional MPC idea from \cite{Jost2015a} to min-max MPC problems resulting in a new robust regional MPC approach. A robust affine feedback law is computed from the solution at the current state and can be reused for subsequent states without solving QPs. We presented two techniques to reduce the number of QPs even further. By this, the online overall computational effort of min-max model predictive control can be reduced significantly. 

\section*{Acknowledgement}
Support by the Deutsche Forschungsgemeinschaft (DFG) under grant MO 1086/15-1 is gratefully acknowledged. We thank Wenjian Wu for his support in checking the results. 

\small
\bibliographystyle{plain}
\bibliography{literature} 

\newpage
\appendix
\section{Numerical data}    

For the min-max and the tube-based approach in Table 1, we used the RPI set $\mathcal{R}=\{x \in \R^n \ | \ T_{\mathcal{R}}x \leq d_{\mathcal{R}} \}$ with 
{\footnotesize
\begin{align*}
T_{\mathcal{R}}=[\text{-}0.334822721594949,	\text{-}0.942281139100084 \\
\text{-}0.317918846692282,	\text{-}0.948117928803084 \\
\text{-}0.299963521841615,	\text{-}0.953950672500615 \\
\text{-}0.253640569521423,	\text{-}0.967298537935858 \\
\text{-}0.229185058677874,\text{-}0.973382868597357 \\
\text{-}0.201195019884297,	\text{-}0.979551205386302 \\
\text{-}0.117871653127393,	\text{-}0.993028838145708 \\
\text{-}0.0663617363880062,	\text{-}0.997795630349005 \\
0.000000000000000,	\text{-}1.000000000000000 \\
0.250417093924644,	\text{-}0.968138047527488 \\
0.440588790533441,	\text{-}0.897709038417393 \\
0.995894097970216,	0.0905259389793321 \\
0.930380886136516,	0.366594335351533 \\
0.751754558173264,	0.659443010627696 \\
0.693428298872540,	0.720525637519398 \\
0.646245605348940,	0.763129489383802 \\
0.567362234739135,	0.823468332476613 \\
0.540329755110452,	0.841453359219796 \\
0.516539177185539,	0.856263556641579 \\
0.471103848089464,	0.882077754121086 \\
0.453378358655099,	0.891318160873663 \\
0.436638430724121,	0.899637082836616 \\
0.401145460614385,	0.916014366387598 \\
0.385834665831679,	0.922567943645646 \\
0.370505236590452,	0.928830377226678 \\
0.334822721594949,	0.942281139100084 \\
0.317918846692282,	0.948117928803084 \\
0.299963521841615,	0.953950672500615 \\
0.253640569521423,	0.967298537935858 \\
0.229185058677874,	0.973382868597357 \\
0.201195019884297,	0.979551205386302 \\
0.117871653127393,	0.993028838145708 \\
0.0663617363880062,	0.997795630349005 \\
0.000000000000000,	1.000000000000000 \\
\text{-}0.250417093924644,	0.968138047527488 \\
\text{-}0.440588790533441,	0.897709038417393 \\
\text{-}0.995894097970216,	\text{-}0.0905259389793321 \\
\text{-}0.930380886136516,	\text{-}0.366594335351533 \\
\text{-}0.751754558173264,	\text{-}0.659443010627696 \\
\text{-}0.693428298872540,	\text{-}0.720525637519398 \\
\text{-}0.646245605348940,	\text{-}0.763129489383802 \\
\text{-}0.567362234739135,	\text{-}0.823468332476613 \\
\text{-}0.540329755110452,	\text{-}0.841453359219796 \\
\text{-}0.516539177185539,	\text{-}0.856263556641579 \\
\text{-}0.471103848089464,	\text{-}0.882077754121086 \\
\text{-}0.453378358655099,	\text{-}0.891318160873663 \\
\text{-}0.436638430724121,	\text{-}0.899637082836616 \\
\text{-}0.401145460614385,	\text{-}0.916014366387598 \\
\text{-}0.385834665831679,	\text{-}0.922567943645646 \\
\text{-}0.370505236590452,	\text{-}0.928830377226678 ],
\end{align*}
}
{\footnotesize
\begin{align*}
d_{\mathcal{R}}=[&0.0693214580803229 \\
&0.0664284061261047\\
&0.0634186283766247\\
&0.0637414504824674\\
&0.0638537776794941\\
&0.0640625898477516\\
&0.0817619105082500\\
&0.0923329894334705\\
&0.105851974067406\\
&0.203197184825767\\
&0.272247311547774\\
&0.394364566746501\\
&0.339042734430747\\
&0.238017010879673\\
&0.208446720595667\\
&0.185219796504134\\
&0.149217644826135\\
&0.137115509713457\\
&0.126574972332281\\
&0.108735935515306\\
&0.101832107009684\\
&0.0953666178493325\\
&0.0846143086314241\\
&0.0799939919533759\\
&0.0754184546935724\\
&0.0693214580803229\\
&0.0664284061261047\\
&0.0634186283766247\\
&0.0637414504824674\\
&0.0638537776794941\\
&0.0640625898477516\\
&0.0817619105082500\\
&0.0923329894334705\\
&0.105851974067406\\
&0.203197184825767\\
&0.272247311547774\\
&0.394364566746501\\
&0.339042734430747\\
&0.238017010879673\\
&0.208446720595667\\
&0.185219796504134\\
&0.149217644826135\\
&0.137115509713457\\
&0.126574972332281\\
&0.108735935515306\\
&0.101832107009684\\
&0.0953666178493325\\
&0.0846143086314241\\
&0.0799939919533759\\
&0.0754184546935724].
\end{align*}
}

\newpage 
Moreover, we used the terminal set $\mathcal{T}=\{x \in \R^n \ | \ T_{\mathcal{T}}x \leq d_{\mathcal{T}} \}$ that has been computed according to \cite{MAYNE2005} such that \begin{align} 
\mathcal{T}=\{x \in \mathcal{S} \mid \forall k \in \N \colon (A+BK_{\infty})^k x \in \mathcal{S}  \}	\nonumber
\end{align}
with $\mathcal{S} \colon = \{ x \in (\mathcal{X} \ominus \mathcal{R}) \mid K_{\infty} x \in (\mathcal{U} \ominus K_{\infty}\mathcal{R}) \}$. The matrices read 
{
\begin{align*}
T_{\mathcal{T}}=[&1.000000000000000,	0.000000000000000 \\
&\text{-}1.000000000000000,	0.000000000000000 \\
&\text{-}0.253574893017553, \text{-}0.967315756943480 \\
&0.253574893017553, 0.967315756943480 \\
&\text{-}0.117740811279278, \text{-}0.993044360217255 \\
&0.117740811279278, 0.993044360217255 ],
\end{align*}
}
\vspace{-10pt}
{
\begin{align*}
d_{\mathcal{T}}=[&9.594531643865977\\
&9.594531643865977\\
&1.170829473660753 \\
&1.170829473660753\\
&2.511330708828057 \\
&2.511330708828057].
\end{align*}
}.

The terminal weighting reads 
{
\begin{align*}
P=[&3.81471424647913,	4.86866526790586 \\
&4.86866526790586,	13.7039014909127]. 
\end{align*}
}

\end{document}